\documentclass{article}
\usepackage{amsfonts,enumerate}

\newtheorem{theorem}{Theorem}

\def\R{{\mathbb R}}
\def\Q{{\mathbb Q}}
\def\L{{\cal L}}
\def\beq{\begin{equation}}
\def\eeq{\end{equation}}
\def\const{{\mathrm{const}}}

\def\rank{{\mathrm{rank}\,}}

\begin{document}

\title{Periodic magnetic geodesics on almost every energy level via
variational methods \thanks{The work was supported by RFBR (grant
09-01-12130-ofi-m) and Max Planck Institute for Mathematics in
Bonn.}}
\author{I.A. Taimanov
\thanks{Institute of Mathematics, 630090 Novosibirsk, Russia;
e-mail: taimanov@math.nsc.ru.} }
\date{}
\maketitle

\medskip

\hfill{\it To V.V. Kozlov on the occasion of his 60th birthday}

\medskip

\section{Introduction}

In the article by using variational methods we prove that for almost
every energy level ``the principle of throwing out cycles'' gives
periodic magnetic geodesics on the critical levels defined by the
``thrown out'' cycles (Theorem \ref{th2}).

In particular, Theorem \ref{th2} implies

\begin{theorem}
\label{th1}
Let $M^n$ be a closed Riemannian manifold and $A$ be a one-form on $M^n$.
If there is a contractible closed curve $\gamma \subset M^n$ such that
\begin{equation}
\label{ineq}
S_E(\gamma) = \int_\gamma \left( \sqrt{E g_{ik} \dot{x}^i\dot{x}^k} +
A_i\dot{x}^i \right) dt < 0,
\end{equation}
then either there is a contractible closed extremal of the functional $S_E$,
either there is a sequence $\gamma_n, n=1,\dots$, of contractible
closed extremals of the functionals $S_{E_n}$ such that
$$
\lim_{n \to \infty} E_n = E, \ \ \
\lim_{n \to \infty} \mathrm{length}(\gamma_n) = \infty.
$$
\end{theorem}

Theorem \ref{th1} is covered by other results obtained by using
methods of symplectic geometry \cite{FS,Contreras} (see \S 3).
However we prove it
by using rather simple variational methods.
Since the formulation of Theorem \ref{th2} needs
some preliminaries we expose it in \S 2.

Now let us first explain what are magnetic geodesics.

Let $M^n$ be a Riemannian manifold and $F$ be a closed two-form on $M^n$.
The motion of a charge in a magnetic field $F$
is given by the Euler--Lagrange equations
$$
\frac{d}{dt}\left( \frac{\partial \L}{\partial \dot{x}^i}\right) -
\frac{\partial \L}{\partial x^i} = 0, \ \ \ i=1,\dots,n,
$$
for the following Lagrangian function
\begin{equation}
\label{lagrange}
\L^\alpha(x,\dot{x}) =
\frac{1}{2}g_{ik}\dot{x}^i\dot{x}^k + A_i^\alpha \dot{x}^i
\end{equation}
where $A^\alpha = d^{-1}F$ are a one-form defined in a domain $U_\alpha$
such that the restriction $F_{U_\alpha}$ of $F$ onto this domain is
an exact form. Therewith the Riemannian metric $g_{ik}$ is used to
evaluate the kinetic energy of the charge.

If $F$ is exact we have a globally defined Lagrangian function. Otherwise
only the Euler--Lagrange equations are globally defined because the
one-form $A$, the vector potential of the magnetic field, enters them via
$F$ and the trajectories of a charge are extremals of
a multivalued functional $W^\alpha = \int \L^\alpha dt$ for which
the variational derivative $\delta W^\alpha$ is globally defined.
The analysis of this picture was done by Novikov
\cite{Novikov81,Novikov82} who started the investigation of variational
problems for such functionals on the spaces of closed curves looking for
periodic trajectories of charges.

The flow corresponding to the Lagrangian (\ref{lagrange})
is called the {\it magnetic geodesic flow.}

The energy
$$
E = \frac{1}{2} g_{ik}\dot{x}^i\dot{x}^k
$$
is preserved along trajectories and on a given energy level the
magnetic geodesic flow is trajectory equivalent to the flow defined by
another Lagrangian function which is homogeneous in velocities:
$$
\L_E^\alpha = \sqrt{E g_{ik}\dot{x}^i \dot{x}^k} + A^\alpha_i\dot{x}^i,
\ \ \ \ E = \mathrm{const}.
$$
The corresponding functionals
$$
S_E^\alpha(\gamma) = \int_\gamma \left(
\sqrt{E g_{ik}\dot{x}^i \dot{x}^k} + A^\alpha_i\dot{x}^i
\right) dt
$$
are defined on the spaces of non-parameterized closed curves on
which they are additive: \beq \label{additive} S(\gamma_1 \cup
\gamma_2) = S(\gamma_1) + S(\gamma_2), \eeq where the curve
$\gamma_1 \cup \gamma_2$ consists in two curves $\gamma_1$ and
$\gamma_2$ which are passed successively.

We thank A. Bahri for helpful discussions.

\section{The principle of throwing out cycles}

The topological problems which arise when we look for some
analogs of the Morse inequalities for
$S_E^\alpha$ were analyzed
by Novikov \cite{Novikov81,Novikov82}. We did not consider here the
case of multivalued functionals and
in the sequel we assume that the magnetic flow is exact, i.e.,
there exists a globally defined one-form $A$ such that
$$
F = dA.
$$

We have a functional \beq \label{functional} S_E(\gamma) =
\int_\gamma \left( \sqrt{E g_{ik} \dot{x}^i\dot{x}^k} + A_i\dot{x}^i
\right) dt \eeq which is defined, as we already mentioned above, on
the space $\Omega^+(M^n)$ of non-parameterized closed curves in
$M^n$. Indeed, since the Lagrangian is homogeneous of degree one in
velocity, the value of $S$ does not depend on the parameterization.

If
$$
\sqrt{E g_{ik} \dot{x}^i\dot{x}^k} +
A_i\dot{x}^i > 0 \ \ \ \mbox{for $\dot{x} \neq 0$},
$$
then the functional $S$ is always positive outside one-point curves and
the number of its extremals can be estimated from below by using the
classical Morse theory. For instance, if all extremals are nondegenerate
(in the Morse sense) we have the simplest version of the Morse
inequalities:
$$
m_k \geq \rank H_k(\Omega^+(M^n),M^n;\Q)
$$
where $m_i$ is the number of extremals of index $k$, $k=0,1,\dots$,
and $M^n \subset \Omega^+(M^n)$ is the set formed by one-point curves.
In general we have the sum of the $k$-th type numbers of extremals
in the left-hand side
(see the recent survey \cite{Taimanov09}).

Assume now that (\ref{ineq}) holds. Then we have to use other
topological reasonings to control critical levels of the functional $S$.
They were introduced by Novikov \cite{Novikov81,Novikov82} as
``the principle of throwing out cycles'' which is as follows:

\begin{itemize}
\item
{\sl Let $w \in H_k(M^n)$ be a nontrivial cycle. We say that it is
``thrown out'' into the set $\{S < 0\}$ if there is a continuous map
$$
F: P \times [0,1] \to \Omega^+(M^n)
$$
such that $P$ is a polyhedron, the map $F_0: P \times 0 \to \Omega^+(M^n)$
is a map into one-point curves which realizes the cycle $w$, i.e.
$F_0(P) \subset M^n$ and $F_{0\,\ast}([P]) = w$, and $F(P \times 1) \subset
\{S < 0\}$. If such a map exists then
it generates a nontrivial cycle in}
$$
H_{k+1}(\Omega^+(M^n),\{S \leq 0\}).
$$
\end{itemize}

For $k=0$ this argument reduces to a version of the mountain-pass lemma,
however it is stronger because it takes care about all cycles and gives rise
to critical levels with $k > 1$.

In \cite{Taimanov83} we proved
that

\begin{itemize}
\item
{\sl if there is a closed contractible
curve $\gamma$ such that $S(\gamma) < 0$ then the whole manifold $M^n$
is thrown out into $\{S <0\}$ which implies the following inequalities:}
\beq
\label{principle}
\rank H_{k+1}(\Omega^+(M^n),M^n) \geq \rank H_k(M^n).
\eeq
\end{itemize}

In fact the construction of such ``throwing out'' works for any functional
$S = \int \L dt$ which satisfy the addition property (\ref{additive}).

If the existence of a critical level implies the existence of a critical point
with the corresponding index on this level, we may put
$m_{k+1}$ into the left-hand side of (\ref{principle}).
However this is not known for the functionals of the type (\ref{functional}).
The usual reasonings by Morse which he did apply for the length functional
\cite{Morse} does not work in this situation (see the discussion in
\cite{Taimanov91}): the gradient deformations for the functional $S$ may
increase the lengths of deformed curves and a priori the gradient deformation
may not converge to critical points.

To every cycle $[z] \in H_k(M^n)$ there corresponds a critical level
$c(w)$ of the functional $S$ which is defined as follows:
\beq
\label{level}
c(w) = \inf_{F} \max_{x \in P \times [0,1]} S(x)
\eeq
where the infimum is taken over all throwings out $F: P \times [0,1] \to
\Omega^+(M^n)$ of the cycle $w$. Moreover for sufficiently small
$\varepsilon >0$ we have
$$
H_{k+1}(\{S \leq c(w) + \varepsilon\},\{S \leq c(w) - \varepsilon\}) \neq 0
$$
which implies (\ref{principle}). These arguments are standard for
the Morse theory.

Now the existence of a critical point on the critical level is
derived as follows: let us take a gradient-type deformation for the
functional $S$ which decreases its values and apply it to the set
$\{S \leq c(w) + \varepsilon\}$. Since it can not be deformed into
$\{S \leq c(w) - \varepsilon\}$ due to topological reasons, this
deformation has either to suspend on a critical point (in a generic
situation of index $k+1$), either there is a subset which realizes a
nontrivial cycle and its deformation diverges to ``infinity'' (in
the finite-dimensional Morse theory on compact manifolds this is
impossible due to compactness however in the non-compact and, in
particular, infinite-dimensional case this results in the existence
of ``a critical point at infinity'' \cite{Bahri}).

If the Palais--Smale
condition \cite{PS}
\footnote{It is said that

{\sl a $C^1$-functional $S$ meets this condition in $\{a \leq S \leq b\}$,
$-\infty \leq a \leq b \leq +\infty$, if for any sequence $x_n$ such that
$|S(x_n)|$ is bounded from above, $S(x_n) \in [2,b]$, and $\lim_{n \to \infty}
|\mathrm{grad}\,S(x_n)|=0$ there is a convergent subsequence.}
}
holds it rules out the latter case and this implies that on every
critical level there is a critical point. Until recently this
condition or some of its reasonable replacements is not established
for the functionals of type (\ref{functional}).

However we show that ``the principle of throwing out cycles''
gives periodic magnetic geodesics on almost every energy level as follows:

\begin{theorem}
\label{th2}
Let $M^n$ be a closed Riemannian manifold and $A$ be a one-form on $M^n$ and
let there is a contractible closed curve $\gamma \subset M^n$ such that
$$
S_E(\gamma) = \int_\gamma \left( \sqrt{E g_{ik} \dot{x}^i\dot{x}^k} +
A_i\dot{x}^i \right) dt < 0.
$$
Then for every nontrivial cycle $w \in H_\ast(M^n)$
either there is a contractible closed extremal $\gamma_0$
of the functional $S_E$ with $S_E(\gamma_0) = c(w)$,
either there is a sequence $\gamma_n, n=1,\dots$, of contractible
closed extremals of the functionals $S_{E_n}$ such that
$$
\lim_{n \to \infty} E_n = E, \ \ \
\lim_{n \to \infty} S_E(\gamma_n) = c(w), \ \ \
\lim_{n \to \infty} \mathrm{length}(\gamma_n) = \infty.
$$
\end{theorem}

Theorem \ref{th1} follows from Theorem \ref{th2} immediately.

\section{On the existence of periodic magnetic geo\-de\-sics}

The study of the periodic problem for magnetic geodesics was initiated
by Novikov \cite{Novikov81,Novikov82}. The principle of throwing out
cycles gives the necessary critical levels in the situation where
the classical Morse theory does not work. However the difficulties with
the Palais--Smale type conditions as we mentioned above hinder a
straightforward derivation of the existence of critical points on these levels.
Until recently the application of the principle of throwing out
cycles to proving the existence of periodic magnetic geodesics
was confirmed in two cases:

1) for magnetic geodesics on a two-torus with everywhere positive magnetic
field $F>0$ \cite{Grinevich}. In these case the magnetic field is exact only on
the universal covering of the torus;

2) for exact magnetic fields which everywhere meet the condition
\beq
\label{bahri}
\min_{|\xi|=1} \left\{\mathrm{Ric}(\xi,\xi) -
\sum_\alpha (\nabla_{e_\alpha} F)(e_\alpha,\xi)\right\} >0
\eeq
where $\mathrm{Ric}$ is the Ricci curvature of the Riemannian metric,
$\xi$ is a tangent vector to $M^n$, and $\{e_\alpha\}$ is the orthonormal basis
in the corresponding tangent space \cite{BT}.
This result is valid for a manifold of any dimension.

There is another result obtained by variational methods in the situation of
strong magnetic fields, when there are throwing out of cycles
\cite{Taimanov84,Taimanov91,Taimanov92}. However
it gives locally minimal periodic magnetic geodesics and reads:

\begin{itemize}
\item
{\sl given a strong magnetic field on a closed oriented two-manifold
there is a non-self-intersecting periodic magnetic geodesic
which is a local minimum of the functional $S$.}
\end{itemize}

An exact field is called strong if there is a two-dimensional
submanifold $N^2$ with boundary such that
$$
\mathrm{length}\,(\partial N^2) +
\int_{N^2} F <0.
$$
There is a similar definition of strong non-exact
magnetic fields \cite{Taimanov91}.

So until recently it was possible either to use the special features of
two-manifolds, either to clarify the Palais--Smale condition under certain
conditions (see (\ref{bahri}) to derive the existence of periodic magnetic
geodesics by variational methods.

There is another approach to proving the existence of such geodesics by
methods of symplectic geometry and it was initiated by Arnold \cite{Arnold}
and Kozlov \cite{Kozlov}. This approach was far-developed and we refer to
the recent article \cite{Ginzburg} which, in particular, contains
an extensive list of references. However we have to stress two results obtained
by symplectic and dynamical methods:

1) in \cite{CMP} it was proved that every exact magnetic field on a closed
two-manifold possesses periodic orbits on all energy levels.
\footnote{Until recently it is the only dimension for which the existence of
a periodic magnetic geodesic on every energy level in an exact magnetic field
is established. In the same dimension there is an example of
a non-exact magnetic field which does not possess periodic magnetic geodesics
on certain (in the example just on one) energy levels \cite{Ginzburg96}.}

The proof splits into three cases:

\begin{enumerate}[(a)]
\item
$E > E_0$ where $E_0$ is some constant, the Man\'e strict critical level.
In this case we have, roughly speaking, the geodesic flow of a Finsler metric
and may apply the classical Morse theory;

\item
$E<E_0$. As it is established in \cite{CMP} this is exactly the case of
strong magnetic fields which is covered by results from
\cite{Taimanov84,Taimanov91};

\item
$E=E_0$. In this case the existence of periodic orbits was established
in \cite{CMP} by methods coming from dynamical systems.
\end{enumerate}

\noindent
Moreover we see that due to \cite{CMP} in the two-dimensional
case the case when
the principle of throwing out cycles gives us necessary critical levels and
the case when the classical Morse theory works are separated exactly by a
single energy level;

2) in \cite{Contreras} by using symplectic and dynamical methods
as well as variational methods applied to the (free time)
action functional which differs from $S$
it was proved for exact magnetic fields
that on an energy level which belong to a total Lebesgue
measure subset of $(0,\infty)$ there is a periodic magnetic geodesic.
This extends the earlier result where that was established for energy
levels belonging to a total measure subset of $(0,d)$ with $d$ is
some small constant \cite{FS}.
\footnote{We have to mention the previous results by Polterovich,
Kerman, and Macarini who established
the existence of periodic magnetic geodesics for sequences
of arbitrary small energy levels provided that the magnetic field is
weakly exact or is given by a symplectic form
(see references in \cite{CMP}.}

Recently two other approaches were introduced: in \cite{Koh} it was proposed
to use the geometric flows for pseudogradient deformations of closed
curves for finding closed magnetic geodesics and a completely new approach
to study of closed magnetic geodesics on surfaces was introduced in
\cite{Schneider}.

The approach used in the present article originates in the comments
from \cite[pp. 192--193]{Taimanov}. For proof we
use the approximation technique from \cite{BT} because the
functional $S$ is even not a $C^1$-functional. Although we discussed
in \cite[pp. 192--193]{Taimanov} that such an approach can be used
for establishing the existence of periodic magnetic geodesics on
energy levels from a total measure subset we did not manage to do
that here by such a mild technique which we use.

\section{Proof of Theorem \ref{th2}}

For brevity we denote
$$
S(\gamma) = S_E(\gamma) = \int_\gamma
\sqrt{E g_{ik}\dot{x}^i \dot{x}^k}\,dt.
$$

Let us consider the family of functionals
$$
S_{\varepsilon,\tau}(\gamma) =
\int_\gamma
\left(
\varepsilon |\dot{x}|^2 + |\dot{x}|^{1+\tau} + A_i\dot{x}^i \right)
dt, \ \ \ \ \varepsilon, \tau \geq 0,
$$
where
$$
\gamma \in H(S^1,M^n),
$$
i.e. $\gamma$ lies in
the Hilbert manifold formed by $H^1$-maps
$$
\gamma: [0,1] \to M^n,
$$
such that $\gamma(0)=\gamma(1)$ \cite{Klingenberg}.

For simplicity, we denote by $L(\gamma)$ the length of $\gamma$:
$$
L(\gamma) = \int_\gamma |\dot{x}|\, dt.
$$

In \cite{BT} it was proved that

\begin{enumerate}[(a)]
\item
the functional $S_{\varepsilon,\tau}$ is a $C^1$-functional on $H(S^1,M^n)$
and it meets the Palais--Smale condition for $\varepsilon,\tau >0$;

\item
the extremals of $S_{\varepsilon,\tau}$ meet the equation
\begin{equation}
\label{deq}
\frac{\partial \gamma^i}{\partial t} + \Gamma^i_{jk}
\dot{\gamma}^j \dot{\gamma}^k =
\frac{g^{ik}F_{kj}\dot{\gamma}^j}{2\varepsilon +
(1+\tau)|\dot{\gamma}|^{\tau-1}},
\end{equation}
or shortly
$$
\frac{D \dot{\gamma}^i}{\partial t} =
\frac{g^{ik}F_{kj}\dot{\gamma}^j}{2\varepsilon +
(1+\tau)|\dot{\gamma}|^{\tau-1}},
$$
and the extremals are arc-length parameterized;

\item
the subset $M^n \subset H(S^1,M^n)$
formed by one-point curves is a manifold of local minima of the functional
$S_{\varepsilon,\tau}$ for $0 \leq \tau <1$;

\item
given a sequence $\gamma_{\varepsilon_n,\tau_n}$ of extremals of
$S_{\varepsilon_n,\tau_n}$ such that their lengths are bounded
by constants $K_0$ and $K_1$:
$$
0 < K_0 \leq L(\gamma_{\varepsilon_n,\tau_n}) \leq K_1 < \infty,
$$
and
$$
\lim_{n \to \infty} (\varepsilon_n,\tau_n) = (\varepsilon_0,\tau_0),
$$
there is a subsequence $\gamma_k$ which converges in $H(S^1,M^n)$
to an extremal $\gamma_\infty$ of the functional $S_{\varepsilon_0,\tau_0}$.
\footnote{In \cite{BT} this was formulated for $(\varepsilon_0,\tau_0) = (0,0)$
however the proof works in the general situation.}
\end{enumerate}

The H\"older inequality implies that
$$
L^m(\gamma) \leq \int_\gamma |\dot{x}|^m\, dt
$$
for any real $m >0$ and the equality is achieved only on arc-length
parameterized
curves: $|\dot{x}| = \const$.
This implies
$$
S_{\varepsilon,\tau} \geq \varepsilon L^2 + L^{1+\tau} + \int A_i\dot{x}^i dt
\geq S.
$$

Let there is a contractible closed curve
$\gamma$ such that
$$
S(\gamma) < 0.
$$
Let $w \in H_k(M^n)$ be a cycle which is thrown out into
$\{S < 0\}$ by the map
$F: P \times [0,1] \to \Omega^+(M^n)$. By continuity this map also defines
the throwings out of $w$ into $\{S_{\varepsilon,\tau} < 0\}$
for sufficiently small
$\varepsilon$ and $\tau$. Hence we have the critical level
$$
c(w) > 0
$$
of $S$ and the corresponding critical levels $c_{\varepsilon,\tau}(w)$
of $S_{\varepsilon,\tau}$. The H\"older inequality implies
$$
c_{\varepsilon,\tau}(w) \geq c_{\varepsilon^\prime,\tau^\prime}(w)
\ \ \ \mbox{for $\varepsilon \geq \varepsilon^\prime, \
\tau \geq \tau^\prime$}.
$$

It is also clear from the continuity that

{\sl for any given $w \in H_\ast(M^n)$ and $\beta > 0$ there exist
$\varepsilon (w,\beta)$ and $\tau (w,\beta)$ such that}
$$
c_{\varepsilon,\tau}(w) \leq c(w) + \beta
\ \ \ \mbox{for $\varepsilon < \varepsilon(w,\beta),
\tau < \tau(w,\beta)$}.
$$

Now let us take a smooth non-decreasing function
$f: \R \to [0,1]$
such that
$$
f(x) = 0 \ \ \mbox{$x \leq \frac{c(w)}{20}$}, \ \
f(x) = 1 \ \ \mbox{$x \geq \frac{c(w)}{10}$}, \ \
0 < f(x) < 1 \ \ \mbox{otherwise},
$$
and instead of $S_{\varepsilon,\tau}$ we consider the $C^1$-functional
$$
F_{w,\varepsilon,\tau} (\gamma) = f(S_{0,\tau}(\gamma)) \cdot
S_{\varepsilon,\tau} (\gamma).
$$

The following statements are clear

\begin{enumerate}[(I)]
\item
for sufficiently small $\varepsilon$ and $\tau$
the throwing out of the cycle $w$ determines a critical level
$$
\widetilde{c}_{\varepsilon,\tau}(w)
$$
which meets the inequality
$$
c_{\varepsilon,\tau}(w) \leq
\widetilde{c}_{\varepsilon,\tau}(w);
$$

\item
for any $\beta>0$ there exist $\varepsilon_1(w,\beta)$
and $\tau_1(w,\beta)$ such that
\beq
\label{3deq}
\widetilde{c}_{\varepsilon,\tau}(w) \leq c(w) + \beta
\ \ \ \mbox{for $\varepsilon < \varepsilon_1(w,\beta),
\tau < \tau_1(w,\beta)$};
\eeq

\item
the functional $F_{w,\varepsilon,\tau}$
meets the Palais--Smale condition in the closed domain
$\{F_{w,\varepsilon,\tau} \geq \frac{c(w)}{10}\}$.
\end{enumerate}

Let us fix some $\beta >0$. The Morse theory
\footnote{Since the considered functionals are only $C^1$,
in this case we have to consider
pseudogradient flows instead of the gradient flows \cite{Rabinowitz}
(see also \cite{Struwe}).}
applied to
the functionals  $F_{w,\varepsilon,\tau}$ implies
for
$$
\varepsilon < \varepsilon_1(w,\beta), \ \ \
\tau < \tau_1(w,\beta)
$$
the existence
of extremals
$\gamma_{\varepsilon,\tau}$
with
$$
F_{w,\varepsilon,\tau}(\gamma_{\varepsilon,\tau}) =
\widetilde{c}_{\varepsilon,\tau},
$$
\beq
\label{bound}
0 \leq \frac{c(w)}{10} \leq L(\gamma_{\varepsilon,\tau}) \leq
\sqrt{\frac{c(w)+\beta}{\varepsilon}}.
\eeq
The latter inequalities follow from (\ref{3deq}) and
the definition of $F_{w,\varepsilon,\tau}$ (in fact, these
functionals were defined exactly for achieving these inequalities).
The curves $\gamma_{\varepsilon,\tau}$
are also extremals of $S_{\varepsilon,\tau}$ and
hence by fixing $\varepsilon>0$ we take a subsequence
$\gamma_{\varepsilon,\tau_n}$ which converges as $\tau_n \to 0$ to a
closed curve $\gamma_\varepsilon$ which the extremal of
$S_{\varepsilon,0}$ meeting the inequality (\ref{bound}). By
(\ref{deq}), these curves $\gamma_\varepsilon$ satisfy the equation
\beq
\label{deq2}
\frac{1}{|\dot{\gamma}_\varepsilon|} \frac{D
\dot{\gamma}_\varepsilon^i}{\partial t} =
\frac{g^{ik}F_{kj}\dot{\gamma}_\varepsilon^j}{2\varepsilon
|\dot{\gamma_\varepsilon}| + 1}.
\eeq
However, since
$\gamma_\varepsilon: [0,1] \to M^n$ is arc-length
parameterized, we have
$$
l(\varepsilon): = |\dot{\gamma}_\varepsilon| = L(\gamma_\varepsilon).
$$
By (\ref{bound}), we have
$$
l(\varepsilon) \leq \sqrt{\frac{c(w)+\beta}{\varepsilon}}.
$$
and therefore
$$
\varepsilon \, l(\varepsilon) \to 0 \ \ \ \mbox{as $\varepsilon \to 0$}.
$$

Now at least one of two cases holds:

1) there is a sequence $\gamma_{\varepsilon_n}$ with uniformly bounded lengths:
$l(\varepsilon_n) \leq K = \const < \infty$. In this case it contains
a subsequence which converges to
a periodic  magnetic geodesic (on the given energy level $E$);

2) there is a sequence $\varepsilon_n$ such that
$$
\nu_n = \varepsilon_n l(\varepsilon_n), \ \
\lim_{n \to \infty}l(\varepsilon_n) = \infty
\ \ \ \mbox{and} \ \ \
\lim_{n \to \infty} \nu_n = 0.
$$
In these case, by (\ref{deq2}),
$\gamma_n = \gamma_{\varepsilon_n}$ are periodic magnetic geodesics
on the energy levels
$E_n = E(1+2\nu_n)$, $\lim_{n \to \infty} S_E(\gamma_n) = c(w)$, and
$\lim_{n \to \infty} \mathrm{length}(\gamma_n) = \infty$.

This proves Theorem \ref{th2}.

\end{document}